# Distributed Robust Optimization Method for AC/MTDC Hybrid Power Systems with DC Network Cognizance

Hai-Xiao Li[1,2], Aleksandra Lekić[1]
[1]*Faculty of Electrical Engineering, Mathematics & Computer Science, Delft University of Technology, Delft, The Netherlands*
[2]*School of Electrical and Electronic Engineering*, *Chongqing University of Technology*, *Chongqing, China*
{H.Li-16, A.Lekic}@tudelft.nl

*Abstract*—AC/multi-terminal DC (MTDC) hybrid power systems have emerged as a solution for the large-scale and long-distance accommodation of power produced by renewable energy systems (RESs). To ensure the optimal operation of such hybrid power systems, this paper addresses three key issues: system operational flexibility, centralized communication limitations, and RES uncertainties. Accordingly, a specific AC/DC optimal power flow (OPF) model and a distributed robust optimization method are proposed. Firstly, we apply a set of linear approximation and convex relaxation techniques to formulate the mixed-integer convex AC/DC OPF model. This model incorporates the DC network-cognizant constraint and enables DC topology reconfiguration. Next, generalized Benders decomposition (GBD) is employed to provide distributed optimization. Enhanced approaches are incorporated into GBD to achieve parallel computation and asynchronous updating. Additionally, the extreme scenario method (ESM) is embedded into the AC/DC OPF model to provide robust decisions to hedge against RES uncertainties. ESM is further extended to align the GBD procedure. Numerical results are finally presented to validate the effectiveness of our proposed method.

*Keywords—Multi-terminal DC grid, DC network cognizance, distributed robust optimization, generalized benders decomposition, extreme scenario method*.

## I. INTRODUCTION

Nowadays, the multi-terminal DC (MTDC) grid composed of voltage-source converters (VSCs) has emerged as a promising solution for transmitting a large amount of power produced by renewable energy systems (RESs) to distant AC grids. Optimal power flow (OPF) is a powerful technique that can benefit VSC-interconnected AC/MTDC hybrid power systems in their economical operation [1]. It is well-known that the fundamental AC/DC OPF problem is nonconvex. The meta-heuristic algorithm [2] and interior point method [3] can be selected to solve such problems. However, the solved results might be stuck in the local optimum. In this regard, convex AC/DC OPF models are developed with a series of linearization and convexification techniques [4]-[6], which can be efficiently solved by off-the-shelf convex solvers. Nevertheless, AC/MTDC hybrid power systems do not favor such a centralized problem-solving approach, due to the substantial communication efforts required for centralized data processing.

Compared to centralized optimization, distributed optimization is preferred. It enables independent decision-making among AC systems (including AC grids and RESs) and DC systems, facilitating coordination through sharing boundary information. Various well-known distributed optimization algorithms were employed in follow-up works to solve AC/DC OPF problems, such as the alternating direction method of multipliers (ADMM) in [7], analytical target cascading (ATC) in [8], and alternating direction inexact newton method (ALADIN) in [9]. However, it is worth noting that the constructed AC/DC OPF problems overlook the issue regarding the DC network switching, which benefits the AC/MTDC hybrid power system in flexible operation [2], [10]. In this case, 0-1 binary variables are inevitably introduced to describe the connection status of DC network lines, resulting in the mixed-integer AC/DC OPF problem. ADMM and ATC lack rigorous convergence behavior for handling such problems, and ALADIN is not scalable to them. Generalized Benders decomposition [11] provides an alternative to solve the mixed-integer AC/DC OPF problem in a distributed manner. Although GBD is archival, its application in hybrid systems still necessitates consideration of some emerging potential issues, such as convergence rate and communication delay.

Additionally, addressing uncertainties from RESs is imperative. A recent popular solution for handling these uncertainties is the two-stage robust optimization based on column and constraint generation (CCG) [12]. Despite CCG involving an iterative solving procedure, it is essentially centralized optimization and requires integration with other distributed optimization algorithms to provide distributed robust solutions, thereby complicating the optimization process. In contrast, the extreme scenario method (ESM)-based two-stage robust optimization is more intuitive and easier to combine with various distributed optimization methods [13], including GBD.

According to the above research review, despite the progress made in OPF modeling, distributed problem-solving, and uncertainty handling, comprehensive studies that address all these aspects together are limited. As previously mentioned, the OPF model for hybrid power systems is expected to be mixed-integer convex when considering DC network cognizance. Hence, it is highly desirable to combine GBD and ESM to achieve distributed robust optimization. Accordingly, our main contributions can be summarized as:
• A mixed-integer convex AC/DC OPF model is explicitly formulated considering the DC network-cognizant constraint.
• GBD is employed to offer distributed optimization, and we enhance the traditional GBD procedure in the aspects of parallel computation and asynchronous updating.
• ESM is employed for robust solutions, and we extend the ESM application to align the GBD procedure.

The paper begins with the mathematical formulation of the AC/DC OPF model that incorporates the DC network-cognizant constraint in Section II. Then, improved GBD is introduced in Section III and ESM with its extended application is presented in Section IV. Numerical studies are presented and discussed in Section V, followed by the conclusion drawn in Section VI.

## II. MIXED-INTEGER CONVEX AC/DC OPF MODEL

In this section, a mixed-integer convex AC/DC OPF model is formulated. This is necessary because the convergence of

This work was supported by CRESYM project Harmony (https://cresym.eu/harmony/).

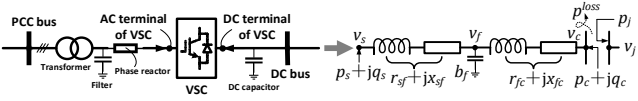

Fig. 1. Equivalent impedance model of the VSC station.

GBD and robustness of ESM cannot be guaranteed if solving nonlinear optimization models.

*A. Linear AC Grid Constraints*

The well-known nonlinear AC power flow can be linearized by the successive linear approximation [14]:

$$p_{ij} = g_{ij}u_i - g_{ij}^{P,k}\frac{u_i+u_j}{2} - b_{ij}^{P,k}(\theta_{ij} - \theta_{ij}^k) + g_{ij}^{P,k}\frac{v_{ij,L}^s}{2}, \quad (1a)$$

$$q_{ij} = b_{ij}^{Q,k}\frac{u_i+u_j}{2} - b_{ij}u_i - g_{ij}^{Q,k}(\theta_{ij} - \theta_{ij}^k) - b_{ij}^{Q,k}\frac{v_{ij,L}^s}{2}, \quad (1b)$$

$$p_i = \sum_{(i,j)} p_{ij} + u_i \sum_j G_{ij}, \quad q_i = \sum_{(i,j)} q_{ij} - u_j \sum_j B_{ij}, \quad (1c)$$

$$\forall i \in \mathbb{m}^{AC}, \quad \forall (i,j) \in \mathbb{I}^{AC}$$

where $\mathbb{m}^{AC}$ and $\mathbb{I}^{AC}$ are used to denote the node and branch sets for the AC grid. $(\cdot)_{ij}$ and $(\cdot)_i$ represent the variable or parameter linked with branch $(i,j)$ and node $i$ of the AC grid, respectively. (1) is the linearized AC power flow equation. $p_{ij}/q_{ij}$ is the real/reactive power flowing on the branch. $p_i/q_i$ is the real/reactive nodal power injection. $u_i \coloneqq v_i^2$ is the squared nodal voltage. $\theta_{ij}$ denotes the phase difference along the branch. $g_{ij}/b_{ij}$ is the conductance/susceptance of the branch. $G_{ij}/B_{ij}$ corresponds to the real/imaginary parts of the nodal admittance matrix. $g_{ij}^{P,k}, b_{ij}^{P,k}, g_{ij}^{Q,k}, b_{ij}^{Q,k}, v_{ij,L}^s, \theta_{ij}^k$ are the determined parameters that are associated with the initial power flow points of the AC grid, whose details can be found in [14, Sec.3].

In addition, some other operational constraints need to be taken into account, including:

$$p_i = p_i^G - p_i^L - p_i^{a2v}, \quad q_i = q_i^G - q_i^L - q_i^{a2v}, \quad (2a)$$

$$0 \le p_i^G \le \overline{p}_i^G, \quad p_i^G \tan(\varphi^{cap}) \le q_i^G \le p_i^G \tan(\varphi^{ind}), \quad (2b)$$

$$-\overline{s}_{ij} \le \cos\left(\frac{n\pi}{2N}\right)p_{ij} + \sin\left(\frac{n\pi}{2N}\right)q_{ij} \le \overline{s}_{ij}, \quad (2c)$$

$$\underline{u}_i \le u_i \le \overline{u}_i, \quad (2d)$$

$$\forall i \in \mathbb{m}^{AC}, \quad \forall (i,j) \in \mathbb{I}^{AC}, \quad n \in \{1, \cdots, N\}$$

where $(\overline{\cdot})/(\underline{\cdot})$ represents the upper/lower bound of variables. (2a) specifies contributions of the AC nodal power injection, including generation production $p_i^G, q_i^G$, load consumption $p_i^L, q_i^L$, and power transmission from the AC grid to the VSC station $p_i^{a2v}, q_i^{a2v}$. (2b) regulates the allowable range for the generator production. $\varphi^{cap}/\varphi^{ind}$ is the capacitive/inductive power factor. (2c) indicates the linearization of $p_{ij}^2 + q_{ij}^2 \le \overline{s}_{ij}^2$, using $n$-polygon approximation [15] to represent the capacity constraint linked with the branch apparent power $s_{ij}$. (2d) indicates the allowable range of $u_i$.

*B. Linear RES Constraints*

RES #$i$ is modeled as an integrated node while the detailed layout is neglected. Hence, we have the below constraints:

$$-\overline{s}_{\#i}^R \le \cos\left(\frac{n\pi}{2N}\right)p_{\#i}^R + \sin\left(\frac{n\pi}{2N}\right)q_{\#i}^R \le \overline{s}_{\#i}^R, \quad (3a)$$

$$0 \le p_{\#i}^R \le \overline{p}_{\#i}^R, \quad p_{\#i}^R = p_{\#i}^{r2v}, q_{\#i}^R = q_{\#i}^{r2v}, \quad (3b)$$

$$\forall \#i \in \{\#1, \cdots, \#I\}, \quad \forall n \in \{1, \cdots, N\}$$

where $(\cdot)_{\#i}$ represents the variable or parameter linked with RES #$i$. $p_{\#i}^R/q_{\#i}^R$ is the real/reactive power produced by the RES. Similar to (2c), (3a) indicates the linearized capacity constraint regarding the RES apparent power output $s_{\#i}^R$. (3b) indicates that $p_{\#i}^R$ is bounded by the maximum available power $\overline{p}_{\#i}^R$, which is affected by natural factors. RES is operated in the grid-connected mode, and the power outputs $p_{\#i}^R, q_{\#i}^R$ equal to the power transmission $p_{\#i}^{r2v}, q_{\#i}^{r2v}$, which flows from the RES to the VSC station.

*C. Mixed-Integer Convex MTDC Grid Constraints*

The well-known nonlinear DC power flow can be convexfied by the second-order cone relaxation [5]:

$$p_i = \sum_{(i,j)} p_{ij}, \quad p_{ij} + p_{ji} = r_{ij}l_{ij}, \quad p_{ij}^2 \le l_{ij}u_i \quad (4a)$$

$$u_i - u_j = r_{ij}(p_{ij} - p_{ji}), \quad (4b)$$

$$i \in \mathbb{m}^{MTDC}, \quad \forall (i,j) \in \mathbb{I}^{MTDC}$$

where $\mathbb{m}^{MTDC}$ and $\mathbb{I}^{MTDC}$ are used to denote the node and branch sets for the MTDC grid. $r_{ij}$ is the resistance of the branch. $l_{ij} \coloneqq i_{ij}^2$ represents the squared branch current. (4) forms the SOC-relaxed power flow constraint by relaxing the instructive equation that $p_{ij} = i_{ij}v_i$.

Some operational constraints need to be considered:

$$\underline{u}_i \le u_i \le \overline{u}_i, \quad p_i = -p_i^{m2v} \quad (5)$$

$$i \in \mathbb{m}^{MTDC}$$

where (5) regulates the allowable range of the squared MTDC nodal voltage $u_i$ and specifies the contribution of MTDC nodal power injection $p_i$, containing the power transmission $p_i^{m2v}$, which flows from the AC grid to the VSC station.

Besides, DC network cognizance is given consideration. The binary variable $\alpha_{ij} \in \{0,1\}$ is introduced to represent the connection status of the DC line. $\alpha_{ij}$ is enabled (binary-1) means the DC line $(i,j)$ is connected and otherwise (binary-0) is disconnected. $\alpha_{ij}$ has the inherent characteristics that $\alpha_{ii} = 0, \alpha_{ij} = \alpha_{ji}$. In this case, (4b) need to be modified, such that:

$$-\mathrm{M}\alpha_{ij} \le p_{ij} \le \mathrm{M}\alpha_{ij}, \quad (6a)$$

$$(u_i - b_{ij}) - (u_j - t_{ij}) = r_{ij}(p_{ij} - p_{ji}), \quad (6b)$$

$$\underline{u}_i(1 - \alpha_{ij}) \le b_{ij} \le \overline{u}_i(1 - \alpha_{ij}), \quad (6c)$$

$$\underline{u}_j(1 - \alpha_{ij}) \le t_{ij} \le \overline{u}_j(1 - \alpha_{ij}), \quad (6d)$$

$$-\underline{u}_i\alpha_{ij} \le u_i - b_{ij} \le \overline{u}_i\alpha_{ij}, \quad -\underline{u}_j\alpha_{ij} \le u_j - t_{ij} \le \overline{u}_j\alpha_{ij}, (6e)$$

$$i \in \mathbb{m}^{MTDC}, \quad \forall (i,j) \in \mathbb{I}^{MTDC}$$

where M represents a large positive number. $b_{ij}, t_{ij}$ are auxiliary variables. It can be deduced that when $\alpha_{ij} = 1$, $u_i - u_j = r_{ij}(p_{ij} - p_{ji})$ is activated and when $\alpha_{ij} = 0$, $u_i - u_j = r_{ij}(p_{ij} - p_{ji})$ is disabled.

*D. Mixed-Integer Convex VSC Constraints*

The nonlinear AC power flow at the AC side of the VSC station ( refers to the PCC bus to the AC terminal, as shown in Fig. 1), can be convexfied by the SOC relaxation[1] proposed in [4]

$$p_i = c_{ii}G_{ii} + \sum_{(i,j)}(c_{ij}G_{ij} - s_{ij}B_{ij}), \quad (7a)$$

$$q_i = -c_{ii}B_{ii} - \sum_{(i,j)}(c_{ij}B_{ij} + s_{ij}G_{ij}), \quad (7b)$$

$$c_{ij} = c_{ji}, s_{ij} = -s_{ji}, c_{ij}^2 + s_{ij}^2 \le c_{ii}c_{jj}, \quad (7c)$$

$$\forall i \in \mathbb{m}^{VSC}, \quad \forall (i,j) \in \mathbb{I}^{VSC}$$

where $\mathbb{m}^{VSC}$ and $\mathbb{I}^{VSC}$ are used to denote the node and branch sets for the AC side of the VSC station. $c_{ii}, c_{ij}, s_{ij}$ are the introduced variables that have links with the squared nodal voltage $u_i$. Given that $u_i \coloneqq e_i^2 + f_i^2$, then we have that $u_i = c_{ii} \coloneqq e_i^2 + f_i^2, c_{ij} \coloneqq e_ie_j + f_if_j, s_{ij} \coloneqq e_if_j - e_jf_i$. SOC relaxation of AC power flow is achieved by relaxing the inst

---

[1] SOC relaxation [4] can be an option for handling the nonlinear power flow of the AC grid. However, commonly used off-the-shelf optimizers sometimes fail to return duals in solving SOC programming (SOCP) problems (https://github.com/jump-dev/Gurobi.jl/issues/217). which affects the application of GBD in the proposed AC/DC OPF model.

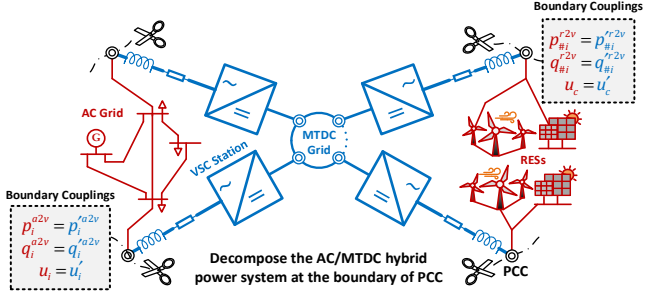

Fig. 2. Illustration regarding the hybrid power system decomposition.

-ructive equation that $c_{ij}^2 + s_{ij}^2 = c_{ii}c_{jj}$.

Some operational constraints at the AC side of the VSC station need to be considered, such that:
$$p_s = p_k^{a2v} \text{ or } p_{\#i}^{r2v}, \quad p_f = 0, \quad p_c = p_j^{m2v} - p_j^{loss}, \quad (8a)$$
$$q_s = q_k^{a2v} \text{ or } q_{\#i}^{r2v}, \quad q_f = u_f b_f, \quad q_c = q_j^{m2v}, \quad (8b)$$
$$\forall k \in \mathbb{m}^{AC}, \quad \{s,f,c\} \in \mathbb{m}^{VSC}, \quad \forall j \in \mathbb{m}^{MTDC}$$
$$\forall \#i \in \{\#1, \cdots, \#I\}$$

where (8) specifies the power injection at the AC side of the VSC station. $p_j^{loss}$ represents the power loss inside the VSC station.

In addition to the aforementioned constraints, a set of constraints are included to characterize the couplings between the AC terminal and DC terminal of the VSC station. The original nonlinear expression is presented below:
$$v_c = \delta v_j, \quad 0 \le \delta \le \overline{\delta}, \quad p_c + p_j + p_j^{loss} = 0, \quad (9a)$$
$$p_j^{loss} = a_{1c}i_c^2 + a_{2c}i_c + a_{3c}, \quad (9b)$$
$$i_c = \sqrt{(p_c^2 + q_c^2)/v_c^2}, \quad 0 \le i_c \le \overline{i}_c, \quad (9c)$$
$$\{s,f,c\} \in \mathbb{m}^{VSC}, \quad \forall j \in \mathbb{m}^{MTDC}$$

where (9a) indicates the voltage and power couplings between the VSC's AC and DC terminals. $\delta$ denotes the PWM modulation factor. (9b) and (9c) indicate that the power loss inside the VSC station is the quadratic function with respect to $i_c$, which is the current at the VSC's AC terminal. $a_{1c}, a_{2c}, a_{3c}$ are the coefficients of the quadratic power loss function. The nonlinear constraints (9) can be handled by piece-wise convex relaxation, and we have that:
$$u_c \le (\overline{\delta})^2 u_j, \quad p_c + p_j + p_j^{loss} = 0, \quad (10a)$$
$$p_j^{loss} = a_{1c}l_c + a_{2c}i_c + a_{3c}, \quad (10b)$$
$$l_c \ge \sum_k i_{c,k}^2, \quad i_c \le \sum_k \{\overline{i}_{c,k}i_{c,k} + \underline{i}_{c,k}i_{c,k} - \overline{i}_{c,k}\underline{i}_{c,k}b_k\}, \quad (10c)$$
$$\sum_k b_k = 1, \sum_k i_{c,k} = i_c, \underline{i}_{c,k}b_k \le i_{c,k} \le \overline{i}_{c,k}b_k, \quad (10d)$$
$$\underline{i}_{c,k} = \frac{\overline{i}_c(k-1)}{K}, \overline{i}_{c,k} = \frac{\overline{i}_c k}{K}, \quad (10e)$$
$$p_c^2 + q_c^2 \le l_c u_c, \quad (10f)$$
$$\{s,f,c\} \in \mathbb{m}^{VSC}, \quad \forall j \in \mathbb{m}^{MTDC}, \quad \forall k \in \{1, \cdots, K\}$$

where the voltage terms $v_c, v_j$ in (9a) is replaced by the squared voltage terms $u_c, u_j$ in (10a) to keep being consistent with the counterparts in (4) and (7). (10a)-(10e) forms a tight quadratic envelope [12] to approximate $l_c \coloneqq i_c^2$ in (9b). $i_{c,k}$ denotes the variable within the subrange $[\underline{i}_{c,k}, \overline{i}_{c,k}]$, and the binary variable $b_k \in \{0,1\}$ is used to denote the status of each subrange, whether it is enabled (binary-1) or disabled (binary-0). (10f) is the SOC relaxation of the equation that $i_c = \sqrt{(p_c^2 + q_c^2)/v_c^2}$ in (9c).

### E. Optimization Objective of the Hybrid Power System

Minimizing the generation costs and the total power losses (including power losses on lines and inside VSC stations) is considered to be the optimization objective for the AC/MTDC hybrid power system, such that:

$$min \left\{ \underbrace{\sum_i \{c_{1i}(p_i^G)^2 + c_{2i}p_i^G + c_{3i}\}}_{generation\ costs} + \underbrace{\sum_i \{p_i^G - p_i^L\} + \sum_i p_{\#i}^R}_{total\ power\ losses} \right\},$$
$$\forall i \in \mathbb{m}^{AC}, \forall \#i \in \{\#1, \cdots, \#I\} \quad (11)$$

where the generation cost is a quadratic function with respect to $p_i^G$. $c_{1i}, c_{2i}, c_{3i}$ are the corresponding coefficients. Total power losses equal the total power generation minus total load demands.

## III. GBD-BASED DISTRIBUTED OPTIMIZATION

Considering the issues regarding the centralized communication burden and data privacy, the AC grid, RESs, and the VSC-MTDC grid should be governed individually. GBD is employed to provide distributed optimization for the coordination among the decomposed sub-systems.

### A. System Decomposition and Multi-Cut Generation

As shown in Fig. 2, the original integrated system can be decomposed into one DC system and several AC systems (refer to the AC grid or RESs). Under this condition, the constructed mixed-integer convex AC/DC OPF model can be reformulated as the below general expression:
$$min\ \{\sum_n \mathcal{F}_{\#n}^{ac}(x_{\#n}, b_{\#n}): h_{\#n}^{cp}(b_{\#n}, b'_{\#n}) = 0,$$
$$x_{\#n} \in X_{\#n}, y \in Y\}, (12)$$
$$\forall \#n \in \{\#1, \cdots, \#N\}$$

where $x_{\#n}$ represents the variables associated with the AC system $\#n$ and $y$ represents the variable associated with the VSC-MTDC grid. $X_{\#n}$ denotes the feasible region of $x_{\#n}$. For the AC system, $X_{\#n}$ is determined by constraints (1)-(3). $Y$ denotes the feasible region of $y$, which is determined by constraints (4)-(10). $b_{\#n}$ represents the boundary variables at PCC. As indicated in (11), the AC system $\#n$ has its own optimization objective $\mathcal{F}_{\#n}^{ac}(x_{\#n}, b_{\#n})$, whereas the VSC-MTDC grid does not. $(\cdot)'$ denote the replicated variables generated after the hybrid power system decomposition, and $h_{\#n}^{cp}(b_{\#n}, b'_{\#n}) = 0$ is the coupling constraints, enforcing consistency among the boundary variables. For instance, the AC grid has the boundary coupling constraints that $p_i^{a2v} = p_i'^{a2v}, q_i^{a2v} = q_i'^{a2v}, u_i = u_i', \forall i \in \mathbb{m}^{AC}$.

Furthermore, (12) can be divided into several slave problems (SPs) related to the AC systems and one master problem (MP) related to the VSC-MTDC grid. During GBD procedure, the original SP $\#n$ is formulated below:
$$OSP \stackrel{\text{def}}{=} min\ \{\mathcal{F}_{\#n}^{ac}(x_{\#n}, b_{\#n}): h_{\#n}^{cp}(b_{\#n}, \hat{b}'^{[m]}_{\#n}) = 0 | \lambda_{\#n},$$
$$x_{\#n} \in X_{\#n}\}, (13)$$

where $(\hat{\cdot})^{[m]}$ denotes the determined variables at the $m$th iteration. $\lambda_{\#n}$ is the dual multiplier corresponding to the constraint that $h_{\#n}^{cp}=0$. (13) provides the upper bound (UB).

However, (13) might have no feasible solutions, and the relaxed one is then formed:
$$RSP \stackrel{\text{def}}{=} min\ \{\|\varepsilon_{\#n}\|_1 + \|\sigma_{\#n}\|_1: h_{\#n}^{cp}(b_{\#n}, \hat{b}'^{[m]}_{\#n}) \le \varepsilon_{\#n} | \mu_{\#n}^\varepsilon,$$
$$-\sigma_{\#n} \le h_{\#n}^{cp}(b_{\#n}, \hat{b}'^{[m]}_{\#n}) | \mu_{\#n}^\sigma, \mu_{\#n}^\varepsilon \ge 0, \mu_{\#n}^\sigma \ge 0,$$
$$x_{\#n} \in X_{\#n}\}, (14)$$

here $\varepsilon_{\#n}$ is the dual multiplier corresponding to $h_{\#n}^{cp} \le \varepsilon_{\#n}$ and $\sigma_{\#n}$ is the dual multiplier corresponding to $-\sigma_{\#n} \le h_{\#n}^{cp}$.

The archival GBD procedure requires that if arbitrary SP $\#n$ has no feasible solution, all SPs should transform into relaxed forms. This necessitates rendering SPs unable to be handled in parallel. To address this issue, we employ the multiple-cut generation [17] in GBD (MGBD). In this case, each SP can generate either one Benders optimality cut originating from (13) or one Benders feasibility cut originating from (14), and multiple cuts are generated in each

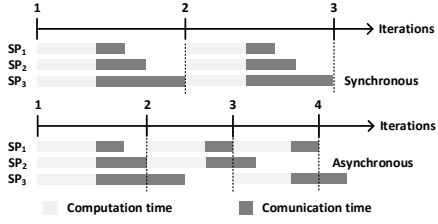

Fig. 3. Illustration regarding the synchronous and asynchronous updates.

---

**Algorithm 1 A**synchronous update of **m**ulti-cut GBD (A-MGBD)

1. Initialize $m \leftarrow 0$, $\hat{b}'^{[0]}_{\#n}, LB^{[0]}_{\#n}, UB^{[0]}_{\#n}$
2. Repeat
   /* SPs Optimization */
3.    **for all** $\#n \in \{\#1,\cdots,\#N\}$ **do**
4.       **if** (13) is feasible **then**
5.          Obtain $\hat{x}^{[m]}_{\#n}, \hat{b}^{[m]}_{\#n}$ by solving (13). Generate Benders optimality cuts to MP, and update $\mathcal{F}^{ac[m+1]}_{\#n}$
6.       **else**
7.          Obtain $\hat{x}^{[m]}_{\#n}, \hat{b}^{[m]}_{\#n}$ by solving (14). Generate Benders feasibility cuts to MP, and keep $\mathcal{F}^{ac[m+1]}_{\#n} \leftarrow \mathcal{F}^{ac[m]}_{\#n}$
8.       **end**
9.    **end**
10.   $UB^{[m+1]} = \sum_n \mathcal{F}^{ac[m+1]}_{\#n}$
    /* MP Optimization */
11.   **while** $\#n \in \mathcal{A}^{[m]}$ finish returning Benders cuts **do**
          Add Benders optimality and feasibility cuts according to (15). Obtain $\hat{b}'^{[m+1]}_{\#n}$ by solving (16) and update $z^{[m+1]}_{\#n}$
12.   **end**
13.   $LB^{[m+1]} = \sum_n z^{[m+1]}_{\#n}$
14.   $m+1 \leftarrow m$
15. **Until** the threshold is met

---

iteration, such that:

$$z_{\#n} \geq \hat{\mathcal{L}}^{[m]}_{\#n} + \nabla \mathcal{L}_{\#n} \cdot (b'_{\#n} - \hat{b}'^{[m]}_{\#n}), \quad (15a)$$
$$0 \geq \hat{\mathcal{H}}^{[m]}_{\#n} + \nabla \mathcal{H}_{\#n} \cdot (b'_{\#n} - \hat{b}'^{[m]}_{\#n}), \quad (15b)$$
$$\mathcal{L}_{\#n} := \mathcal{F}^{ac}_{\#n}(x_{\#n}, b_{\#n}) + \lambda^\top_{\#n} h^{cp}_{\#n}(b_{\#n}, b'_{\#n}) \quad (15c)$$
$$\mathcal{H}_{\#n} := \mu^{\varepsilon\top}_{\#n} \hbar_{\#n}(b_{\#n}, b'_{\#n}) - \mu^{\sigma\top}_{\#n} h^{cp}_{\#n}(b_{\#n}, b'_{\#n}) \quad (15d)$$

Consequently, MP is formulated with the returened Benders cuts, providing the lower bound (LB), such that:

$$MP \overset{\text{def}}{=} \min\{\sum_n z_{\#n} : Eqs.(15a) \text{ and } (15b)\}, \quad (16)$$
$$\forall \#n \in \{\#1,\cdots,\#N\}, \quad \forall [m] \in \{[1],\cdots,[M]\}$$

Eventually, the GBD procedure is implemented by iteratively solving SPs and MP. The update process is terminated until the residuals of coupling constraints for SPs and MP are smaller than the threshold.

### B. Asynchronous Updating

We further consider the communication delay during the GBD procedure. As shown in Fig. 3, affected by different communication delays, MP is unable to receive all returned Benders cuts simultaneously. Synchronous and asynchronous updating can be selected to respond to this situation. The synchronous updating requires waiting for all SPs to return Benders cuts. In this way, the total computation time is limited by the "slowest" SP-solving. In contrast, the asynchronous updating only needs a minimum number of $n \geq 1$ SPs to finish cut-returning. To ensure sufficient freshness, every SP must finish updating at least once every $m \geq 1$ iteration. For the asynchronous updating steps of MGBD, we denote the set of SPs that finish returning Benders cuts at the $m$th iteration as $\mathcal{A}^{[m]}$ and remainings as $\mathcal{C}^{[m]}$. Hence, we have that $\mathcal{A}^{[m]} \cap \mathcal{C}^{[m]} = \emptyset, \mathcal{A}^{[m]} \cup \mathcal{C}^{[m]} = \{\#1,\cdots,\#N\}$. The update steps are described in **Algorithm 1**.

## IV. ESM-BASED ROBUST SOLUTION

To address uncertainties from RESs, the aforementioned deterministic OPF (DOPF) model (formulated by (1)-(11) is expanded into the ESM-based two-stage robust optimization model (ROPF). Furthermore, this ESM-based ROPF model is further extended to align the GBD procedure.

### A. Two-Stage Robust Decision-Making With ESM

The DOPF model is expanded into the two-stage ROPF model, such that:

$$\min\left\{\frac{1}{|\mathbb{S}|}\sum_s \mathcal{F}^{hyb}_{\{s\}}(\xi_{\{s\}}, v, w_{\{s\}}): \hbar^{hyb}_{\{s\}}(\xi_{\{s\}}, v, w_{\{s\}}) = 0,\right.$$
$$\left. \mathcal{g}^{hyb}_{\{s\}}(\xi_{\{s\}}, v, w_{\{s\}}) \leq 0\right\} \quad (17)$$
$$\forall s \in \mathbb{S}$$

where $(\cdot)_{\{s\}}$ denotes the variables linked with the specific possible scenario $s$. $\mathcal{F}^{hyb}_{\{s\}}, \hbar^{hyb}_{\{s\}}, \mathcal{g}^{hyb}_{\{s\}}$ respectively represent the optimization objective, equality constraint, and inequality constraint of the hybrid power system in the possible scenario $s$. $\xi_{\{s\}}$ means the uncertain variable in the scenario $s$, referring to the maximum RES available power $\overline{p}^R_{\#i\{s\}}$. $v$ refer to the "here-and-now" variables determined in the first stage before $\xi_{\{s\}}$ is revealed, including DC network connection status $\alpha_{ij}$ and generator power outputs $p^G_i, q^G_i$. The remaning variables are "wait-to-see" variables determined in the second stage after $\xi_{\{s\}}$ is revealed. They are uniformly denoted as $w_{\{s\}}$. As indicated in (17), the robust $y$ can be obtained if we solve (17) by adding constraints involving the collection of all possible scenarios $\mathbb{S}$, but it is evidently computationally massive.

ESM provides an alternative approach to solve (17) by considering the collection of all extreme scenarios $\mathbb{e}$ instead of $\mathbb{S}$. Accordingly, (10) is reduced to:

$$\min\left\{\frac{1}{|\mathbb{e}|}\sum_e \mathcal{F}^{hyb}_{\{e\}}(\xi_{\{e\}}, v, w_{\{e\}}): \hbar^{hyb}_{\{e\}}(\xi_{\{e\}}, v, w_{\{e\}}) = 0,\right.$$
$$\left. \mathcal{g}^{hyb}_{\{e\}}(\xi_{\{e\}}, v, w_{\{e\}}) \leq 0\right\} \quad (18)$$
$$\forall e \in \mathbb{e} \subseteq \mathbb{S}$$

where $(\cdot)_{\{e\}}$ denotes the variables linked with the specific extreme scenario $e$. We can konow that solving (18) requires much less computaitnal effort compared to solving (10) due to $|\mathbb{e}| \ll |\mathbb{S}|$. The effectivness of ESM holds rquiring that $\hbar^{hyb}$ is linear and $\mathcal{g}^{hyb}$ is linear or SOC relaxed. The proof can be found in [13, Sec.4]. As described in (1)-(10), the proposed AC/DC OPF model meets this requirement.

### B. Extension of ROPF for Distributed Problem-Solving

With ESM, the original SP #$n$ in (13) is modified by considering all extreme scenarios $\mathbb{e}$, such that:

$$OSP \overset{\text{def}}{=} \min\left\{\frac{1}{|\mathbb{e}|}\sum_e \mathcal{F}^{ac}_{\#n}\left(\xi^{ac}_{\#n\{e\}}, v^{ac}_{\#n}, \widetilde{w}^{ac}_{\#n\{e\}}, b_{\#n\{e\}}\right):\right.$$
$$h^{cp}_{\#n\{e\}}\left(b_{\#n\{e\}}, \hat{b}'^{[m]}_{\#n\{e\}}\right) = 0, \hbar^{ac}_{\#n\{e\}}(\xi^{ac}_{\#n\{e\}}, v^{ac}_{\#n}, \widetilde{w}^{ac}_{\#n\{e\}}) = 0,$$
$$\left.\mathcal{g}^{ac}_{\#n\{e\}}(\xi^{ac}_{\#n\{e\}}, v^{ac}_{\#n}, \widetilde{w}^{ac}_{\#n\{e\}}) \leq 0\right\}, (19)$$
$$\forall \#n \in \{\#1,\cdots,\#N\}, \quad \forall e \in \mathbb{e}$$

where $\xi_{\{e\}}, v, w_{\{e\}}$ are extended to $\xi^{ac}_{\#n\{e\}}, v^{ac}_{\#n}, w^{ac}_{\#n\{e\}} := [\widetilde{w}^{ac}_{\#n\{e\}}, b_{\#n\{e\}}]$, respectively representing the uncertain variables, first-stage varaibles, and second-stage variables related to the AC system, in the extreme scenario $e$. Particularly, $\widetilde{w}^{ac}_{\#n\{e\}}$ refers to the second-stage vairables $w^{ac}_{\#n\{e\}}$ exclude the boundary variables $b_{\#n\{e\}}$. (14) and (16) can also be modified by this kind of variable extension.

However, from the perspective of distributed problem-solving, GBD solving (19) needs extra communication. We

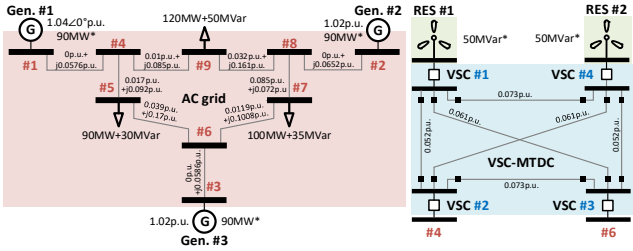

Fig. 4. Test system of the AC/MTDC hybrid power system.

TABLE I. COMPARISON BETWEEN DOPF ROPF AND E-ROPF MODELS IN HANDLING POSSIBLE SCENARIOS

| Optimization model | Optimization objective in 100 possible scenarios | | | |
|---|---|---|---|---|
| | Feasible ratio | Maximum | Minimum | Average |
| DOPF | 73% | - | - | - |
| ROPF | 100% | 8.3114 | 8.3109 | 8.3110 |
| E-ROPF | 100% | 8.3117 | 8.3117 | 8.3117 |

✧ Models are solved via centralized optimization. The feasible ratio means the proportion of possible scenarios that have feasible solutions.

TABLE II. COMPARISON OF E-ROPF WITH AND WITHOUT DC NETWORK COGNIZANCE

| Network-cognizant constraint | DC topology | VSC connection | Nodal voltage constraint |
|---|---|---|---|
| Without | Fixed | (1,2), (2,3), (3,4) (4,1) | Violate |
| With | Flexible | (1,2), (2,3), (3,4), (4,1), **(1,3), (2,4)** | Satisfy |

✧ Models are solved via centralized optimization.

TABLE III. COMPARISON OF GBD AND MGBD IN SOLVING E-ROPF

| Approaches | Iterations | $\|UB-LB\|/\|LB\| \times 100\%$ | $\sigma \times 100\%$ |
|---|---|---|---|
| GBD | 23 | 0.0240% | 0.0172% |
| MGBD | 15 | 1.4584% | 0.0181% |

✧ $\sigma$ refers to the relative error between the solved result by distributed approaches and the benchmark, which is obtained by solving the E-ROPF via centralized optimization.

denote $\mathbb{e}_{\#n}$ as the local extreme scenarios induced by the uncertain variable $\xi^{ac}_{\#n}$ related to AC system $\#n$. In this way, we have that $|\mathbb{e}| = \prod_n |\mathbb{e}_{\#n}|$ where $\mathbb{e}$ can be interpreted as the arrangement of $\mathbb{e}_{\#n}$. It indicates that besides the local $\mathbb{e}_{\#n}$, AC system $\#n$ needs to be informed the remaining extreme scenarios $\mathbb{e}_{\#r}, r \neq n$. Considering this issue, a consensus condition is designed, i.e., the boundary variables of each AC system need to remain consistent across all extreme scenarios $\mathbb{e}$. Consequently, (19) is transformed, such that:

$$OSP \overset{\text{def}}{=} \min \left\{ \frac{1}{|\mathbb{e}_{\#n}|} \sum_e \mathcal{F}^{ac}_{\#n}\left(\xi^{ac}_{\#n\{e\}}, v^{ac}_{\#n}, b_{\#n}, \widetilde{w}^{ac}_{\#n\{e\}}\right): \right.$$
$$h^{cp}_{\#n}\left(b_{\#n}, \hat{b}'^{[m]}_{\#n}\right) = 0, h^{ac}_{\#n\{e\}}\left(\xi^{ac}_{\#n\{e\}}, v^{ac}_{\#n}, \widetilde{w}^{ac}_{\#n\{e\}}\right) = 0,$$
$$\left. g^{ac}_{\#n\{e\}}\left(\xi^{ac}_{\#n\{e\}}, v^{ac}_{\#n}, \widetilde{w}^{ac}_{\#n\{e\}}\right) \leq 0 \right\}, (20)$$
$$\forall \#n \in \{\#1, \cdots, \#N\}, \qquad \forall e \in \mathbb{e}_{\#n}$$

where the boundary variables $b_{\#n\{e\}}$, as the second-stage variables, are coverted into $b_{\#n}$, as the first-stage variables. $b_{\#n\{e\}}$ in (14) and (16) can also be converted in the similary way. Such distributed problem-solving is equivalent to adding an additional compulsory constraint that $b_{\#n\{e\}} = b_{\#n\{e^*\}} \cdots = b_{\#n\{e^{*\cdots*}\}}, \forall e, e^*, e^{*\cdots*} \in \mathbb{e}$ in (18), forming an extended ROPF (E-ROPF) model.

Comparing (19) and (20), we can know that the E-ROPF model enables SP $\#n$ to consider the collection of local extreme scenarios $\mathbb{e}_{\#n}$ instead of the collection of all extreme scenarios $\mathbb{e}$, avoiding extra communication costs. Details regarding the E-ROPF model can be seen in [13, Sec.4].

## V. CASE STUDY

We use the test system shown in Fig. 4, and the based power and voltage are respectively 100MW and 345kV. The system-wide power and voltage are in the calculation of per-unit values. We set the tighter voltage bounds as $\underline{u}_i := 0.955 p.u.$, $\overline{u}_i := 1.045 p.u., \forall i \in \mathbb{m}^{AC} \cup \mathbb{m}^{VSC} \cup \mathbb{m}^{MTDC}$, for preventing the voltage violations indused by power flow errors. Other key bounds are set as $\arccos(\varphi^{cap}) = \arccos(\varphi^{ind}) := 0.9$, $\overline{s}_{ij} := 1 p.u., (i,j) \in \mathbb{l}^{AC}, \overline{i}_c := 1 p.u., c \in \mathbb{m}^{VSC}, \overline{\delta} := 1$. For GBD, MGBD, A-MGBD, the residual threshold is set to $1 \times 10^{-5}$. The case study is coded on the MATLAB platform. YALMIP toolbox is utilized to provide a mathematic modeling environment. CPLEX is selected as the solver.

### A. Validation of ESM-Based ROPF Model

We assume that the uncertain variables $\bar{p}^R_{\#1}$ and $\bar{p}^R_{\#2}$ are respectively within the interval that $[0.3, 0.5] p.u.$ and $[0.2, 0.5] p.u.$. 100 possible scenarios related to $\bar{p}^R_{\#1}, \bar{p}^R_{\#2}$ are randomly generated. The conventional DOPF model considers one deterministic scenario that is $(\bar{p}^R_{\#1}, \bar{p}^R_{\#2}) := (0.5, 0.5) p.u.$. In contrast, the ESM-based ROPF model considers four extreme scenarios that are $(\bar{p}^R_{\#1}, \bar{p}^R_{\#2}) \in \{(0.5, 0.5), (0.5, 0.2), (0.3, 0.2), (0.3, 0.2)\} p.u.$. E-ROPF model considers two sets of local extreme scenaiors that are $\bar{p}^R_{\#1} \in \{0.5, 0.3\} p.u., \bar{p}^R_{\#2} \in \{0.5, 0.2\} p.u.$. As presented in Table I, the decision-making based on DOPF model is not suitable for all generated possible scenarios, and it will result in infeasible solutions in 27% scenarios. In contrast, the decision-making based on both ROPF and EROPF models are feasible for all generated possible scenarios by only considering extreme scenarios, and the computation cost is cheap. Taking the determined $\hat{\alpha}_{ij}, \forall(i,j) \in \mathbb{l}^{MTDC}$ and $\hat{p}^G_i, \hat{q}^G_i$ $\forall i \in \mathbb{m}^{AC}$ as constants into the generated possible scenarios to validate the performances of ROPF and E-ROPF models on hedging against uncertainties (the E-ROPF model additionally take the determined boundary variables $b_{\#n}$ as constants into the generated possible scenarios). We find that the E-ROPF model provides the slightly more conservative results than the ROPF model. This is because the E-ROPF model has one additional compulsory constraint. However, the E-ROPF model guarantees feasibility for the arbitrary possible scenarios and benefits reducing the communication burden when using GBD.

### B. Influence of the DC Network-Cognizant Constraint

We investigate the influence of the DC network-cognizant constraint formulated in (6) on the E-ROPF model. As presented in Table II, if the DC network cognizance is not considered, it implies that the DC network topology is fixed. In this case, the E-ROPF model has no feasible solution since the nodal voltage constraint is violated. When considering network-cognizant constraint, the E-ROPF model enables flexibly reconfiguring DC topology to satisfy the nodal voltage constraint, and feasible solutions can be found.

### C. Discussion Regarding Multi-Cut Generation and Asynchronous Updating in GBD

We first compare GBD and MGBD in solving the E-ROPF model to investigate the effect of multi-cut generation in GBD. As shown in Table III, MGBD converges faster than GBD, thanks to MGBD returning more Benders cuts in each iteration. Regarding the optimality, for the given case study, it can be considered that GBD and MGBD exhibit almost eq

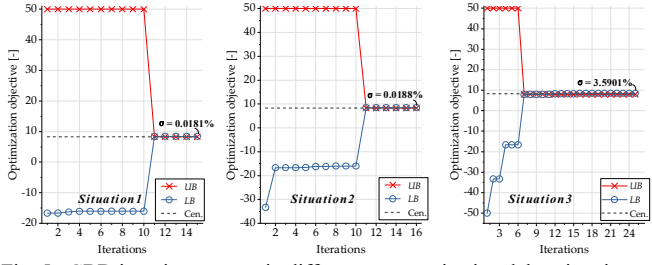

Fig. 5. GBD iteration process in different communication delay situations.

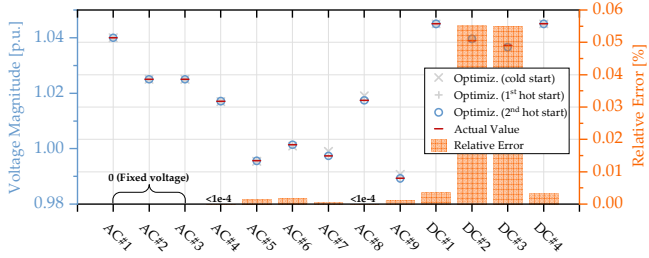

Fig. 6. The calculated system-wide voltage profiles (based on the deterministic scenario and the fixed DC topology). The actual values are computed by using IPOPT to solve the original nonconvex AC/DC OPF model with dertermined $\hat{p}_i^G, \hat{q}_i^G, \hat{p}_{\#1}^R, \hat{q}_{\#1}^R, \hat{p}_{\#2}^R, \hat{q}_{\#2}^R, \hat{\delta}$ which are obtained by solving the proposed mixed-integer convex AC/DC OPF model.

-uivalent performance as both $\sigma$ are pretty close.

We further discuss the asynchronous updating in GBD. Particularly, we set three specific communication delay situations. **Situation1**: In every iteration, all SPs must have been solved and returned cuts to MP, and then the iteration is activated. This situation can be regarded as synchronous updating. **Situation2**: In every iteration, at least two SPs have already been solved, and then the iteration is activated. The time cost for solving SPs associated with the AC grid, RES#1, RES#2 are assumed to be 1:1:2. This situation belongs to asynchronous updating. **Situation3**: In every iteration, at least two SPs have been solved, and then the iteration is activated. The time cost for solving SPs associated with the AC grid, RES#1, RES#2 are assumed to be 1:2:4. This situation also belongs to the asynchronous updating. As presented in Fig. 5, asynchronous updating still converges successfully, and compared with synchronous updating, it needs more iterations but saves time in every iteration. Hence, in terms of the whole iteration process, the total time consumed by A-MGBD might be less. However, we can observe that the converged result in **Situation3** exhibits an obvious deviation. It implies that as communication delays increase, the optimality of A-MGBD is substantially affected.

*D. Test Regarding Power Flow Accuracy*

Our work employs a series of approximations to handle the nonlinear power flow constraints, and it is necessary to test the power flow accuracy after approximations. The initial power flow points in (1) are obtained with the base power flow status and then changed to the updated power flow status (power injection at AC nodes #4 and #6 are changed from zeros to the optimized value. Generation power outputs are changed from the base values to the optimized value). As shown in Fig. 6, after three initial point updates, the proposed mixed-integer convex AC/DC OPF model exhibits high power flow accuracy. However, it is a trade-off between power flow accuracy and computation burden induced by the initial power flow point updates.

## VI. Conclusion

In this paper, a mixed-integer convex AC/DC OPF model considering DC network cognizance is constructed for the AC/MTDC hybrid power system. The system-wide flexible operation is achieved through optimal DC topology reconfiguration. The improved GBD and extended ESM are combined to achieve distributed robust decision-making when confronted with uncertain scenarios related to RESs. However, this current work does not consider the impacts of VSC local control on the hybrid power system. In future work, we will further expand the OPF model, embedding the VSC droop control constraint.